\documentclass[a4paper,12pt]{article}
\usepackage[dvips]{epsfig}
\usepackage{amsmath,amssymb,euscript,graphicx,amsfonts}
\usepackage{enumerate,color}
\usepackage{graphicx}
\usepackage{hyperref}
\usepackage{tikz}

\newtheorem{theorem}{Theorem}[section]
\newtheorem{proposition}[theorem]{Proposition}
\newtheorem{lemma}[theorem]{Lemma}
\newtheorem{corollary}[theorem]{Corollary}

\newcommand{\Qed}{\rule{2.5mm}{3mm}}

\newcommand{\B}{\mathcal{B}}

\newcommand{\ZZ}{\mathbb{Z}}

\newcommand{\FF}{\mathbb{F}}

\newtheorem{hypothesis}[theorem]{Hypothesis}%
\newcounter{case}

\renewcommand{\thecase}{\arabic{case}}

\newcounter{subcase}

\numberwithin{subcase}{case}

\def\soc{\hbox{\rm Soc}}

\def\PGL{\hbox{\rm PGL}}
\def\GL{\hbox{\rm GL}}

\def\SL{\hbox{\rm SL}}

\def\PSL{\hbox{\rm PSL}}

\def\demo{{\bf Proof}\hskip10pt}
\def\di{\bigm|} \def\lg{\langle} \def\rg{\rangle}
\def\qqed{\hfill $\Box$}
\def\a{\alpha} \def\b{\beta}

\def\Aut{\hbox{\rm Aut\,}}
  
\def\f{\noindent}

\def\a{\alpha} \def\b{\beta}

 \def\og{\overline G} \def\oh{\overline H}

   \def\o1{\overline 1}

 \def\olx{\overline x} 

\def\o{\overline}   
\def\di{\bigm|} \def\lg{\langle} \def\rg{\rangle}

\def\di{\bigm|} \def\lg{\langle} \def\rg{\rangle}

\usepackage{minitoc}

\begin{document}

	\begin{center}
		{\bf\large  On Hamilton cycles in quasiprimitive graphs of order 
			$2pq$}
		\footnote{This work was supported in part  by the National Natural Science Foundation of China
			(12471332).}
	\end{center}
	
	
	\begin{center}
		Hao Yu$^{1}$ and Tianlei Zhou$^{2,*}$\\
		\medskip
		
		\begin{minipage}[t]{0.47\textwidth}
			\centering
			{\it\small
				$^{1}$School of Mathematics, Guangxi University,\\
				Guangxi Base, Tianyuan Mathematical Center in Southwest China,\\
				Guangxi Center for Mathematical Research,\\
				Nanning 530004, P. R. China
			}
		\end{minipage}
		\hfill
		\begin{minipage}[t]{0.47\textwidth}
			\centering
			{\it\small
				$^{2}$School of Mathematical Sciences,\\
				Capital Normal University,\\
				Beijing 100048,\\
				P. R. China
			}
		\end{minipage}
	\end{center}
	
	\begingroup
	\renewcommand{\thefootnote}{*}
	\footnotetext{Corresponding author: 2230501007@cnu.edu.cn}
	\endgroup
	
	\renewcommand{\thefootnote}{\empty}

\begin{abstract}
	Let $p>q>5$ be primes. We determine the socles and corresponding
	actions of quasiprimitive imprimitive groups of degree $2pq$,
	all of which are almost simple. Using this classification,
	we prove that a connected graph of order $2pq$ admitting a
	quasiprimitive automorphism group contains a Hamilton cycle,
	except possibly for three specified families of actions with
	socle $\PSL(2,p)$ or $\PSL(2,q^2)$.
\end{abstract}

\section{Introduction}
\label{sec:intro}
\indent

In 1969, Lov\'asz \cite{L70} asked whether every finite connected
vertex-transitive graph contains a Hamilton path. In 1981,
Alspach \cite{A83} asked whether there are infinitely many
connected vertex-transitive graphs without a Hamilton cycle.
No connected vertex-transitive graph without a Hamilton path is known.
Only four connected vertex-transitive graphs on at least three
vertices are known to have no Hamilton cycle: the Petersen graph,
the Coxeter graph, and their truncations, obtained by replacing
each vertex with a triangle. None of these four graphs is a
Cayley graph. A related conjecture asserts that every connected
Cayley graph on at least three vertices contains a Hamilton cycle;
see \cite{CG96,D83,GWM14,GM07,GKMM12,KM09,MR234,DM83,WM18,WM15}.

Hamilton paths and, in many cases, Hamilton cycles have been
shown to exist in various families of vertex-transitive
graphs of orders $kp$, $p^j$, and $2p^2$, where $p$ is prime;
see \cite{C98,MR4328721,MR4548744,KM08,KMZ12,KS09,DM87,DM92,DM85,MP82,MP83,MS4,MS2,MS5,Z15}
and the survey \cite{KM09}. In 2021, Du, Kutnar and
Maru\v{s}i\v{c} \cite{MR4328721} resolved the Hamilton cycle
problem for connected vertex-transitive graphs of order $pq$,
where $p$ and $q$ are distinct primes: every such graph
contains a Hamilton cycle, except for the Petersen graph.

For graphs of order $2pq$, Du, Tian and Yu \cite{MR4548744}
proved that every connected vertex-transitive graph whose
automorphism group acts primitively on vertices contains a
Hamilton cycle, except for the Coxeter graph. Du and Zhou
\cite{6p} proved that every connected vertex-transitive graph
of order $6p$, where $p$ is prime, contains a Hamilton cycle,
except for the truncation of the Petersen graph. Zhou \cite{Zh}
obtained Hamiltonicity results for non-quasiprimitive graphs
of order $2pq$ under hypotheses involving maximal intransitive
normal subgroups with orbits of prime length.

We consider the quasiprimitive case for primes $p>q>5$.
Every nontrivial normal subgroup in a quasiprimitive action
is transitive, so reduction by intransitive normal subgroups
is unavailable. With the primitive case covered by the result
above, we focus on imprimitive actions.
Our first result classifies quasiprimitive imprimitive groups
of degree $2pq$ in terms of their socles and actions.
All such groups are almost simple. The same socle can occur
with different stabilizers, and the resulting actions need not
be equivalent. Keeping these actions distinct is necessary for
the subsequent graph analysis. We determine their possible
socles and the point stabilizers of the restricted actions,
using the classifications of primitive groups of prime degree
or degree a product of two primes and the subgroups of their
point stabilizers. The resulting list also supplies the group
actions needed for the graph problem.

Our second result uses this classification to prove Hamiltonicity
for connected quasiprimitive graphs of order $2pq$, except
possibly for three specified families of actions with socle
$\PSL(2,p)$ or $\PSL(2,q^2)$. The proofs use coset graph
descriptions of these actions. We state the two results below.

\begin{theorem}
	\label{th1.1}
	Let $p$ and $q$ be primes with $p>q>5$, and let $G$ be a
	quasiprimitive imprimitive permutation group of degree $2pq$.
	Then $G$ is almost simple, and its socle and the corresponding
	restricted action are as listed in Table~\ref{tab5}.
\end{theorem}

\begin{theorem}
	\label{the:main}
	Let $\Gamma$ be a connected graph of order $2pq$, where
	$p>q>5$ are primes, and let $G\leq\Aut(\Gamma)$ act
	quasiprimitively and imprimitively on $V(\Gamma)$. Then $\Gamma$ contains
	a Hamilton cycle, except possibly when the action of $G$
	on $V(\Gamma)$ is one of those listed in Rows~$6$--$8$
	of Table~\ref{tab5}.
\end{theorem}

\begin{table}[htbp]
	\centering
	\caption{Socles and restricted actions for quasiprimitive
		imprimitive groups of degree $2pq$, where $5<q<p$ are primes}
	\label{tab5}
	\renewcommand{\arraystretch}{1.18}
	\resizebox{\textwidth}{!}{%
		\begin{tabular}{|c|l|c|l|l|}
			\hline
			Row & $T=\operatorname{soc}(G)$ & $(p,q)$ or degree condition
			& Action of $T$ & Conditions or remarks\\
			\hline
			1 & $A_p$ & $(p,(p-1)/2)$
			& cosets of $A_{p-2}$
			& ordered pairs of distinct points of a $p$-set\\
			\hline
			2 & $A_{p+1}$ & $(p,(p+1)/2)$
			& cosets of $A_{p-1}$
			& ordered pairs of distinct points of a $(p+1)$-set\\
			\hline
			3 & $\PSL(2,r^n)$
			& $((r^n+1)/2,q)$
			& cosets of $\ZZ_r^{\,n}\rtimes \ZZ_{(r^n-1)/(2q)}$
			& $r$ an odd prime; $q\mid r^n-1$\\
			\hline
			4 & $\PSL(d,r^n)$
			& $pq=\dfrac{r^{nd}-1}{r^n-1}$
			& cosets of $\ZZ_r^{\,n(d-1)}\rtimes L$
			& $d\geq3$; $[\GL(d-1,r^n)/Z:L]=2$\\
			\hline
			5 & $\PSL(d,r^n)$
			& $p=\dfrac{r^{nd}-1}{r^n-1}$
			& cosets of $\ZZ_r^{\,n(d-1)}\rtimes L$
			& $d\geq3$; $[\GL(d-1,r^n):L]=2q$\\
			\hline
			6 & $\PSL(2,r^2)$
			& $((r^2+1)/2,r)$
			& cosets of $\PSL(2,r)$
			& $r$ an odd prime\\
			\hline
			7 & $\PSL(2,p)$
			& $(p,(p-\varepsilon)/2)$
			& cosets of $D_{(p+\varepsilon)/2}$
			& $\varepsilon\in\{1,-1\}$\\
			\hline
			8 & $\PSL(2,p)$
			& $(p,(p-\varepsilon)/2)$
			& cosets of $\ZZ_{(p+\varepsilon)/2}$
			& $\varepsilon\in\{1,-1\}$\\
			\hline
			9 & $\PSL(2,23)$ & $(23,11)$
			& cosets of $A_4$ & $A_4<S_4<T$\\
			\hline
			10 & $M_{23}$ & $(23,11)$
			& cosets of $\PSL(3,4)$
			& ordered pairs of distinct points in its natural action\\
			\hline
	\end{tabular}}
\end{table}
Rows~$6$--$8$ of Table~\ref{tab5} are the three possible exceptional
families in Theorem~\ref{the:main}. They correspond, respectively,
to actions on cosets of $\PSL(2,q)$, a dihedral subgroup, and a
cyclic subgroup.

\section{Terminology, notation and some useful results}
\label{sec:pre}
This section consists of four subsections dealing with basic definitions and notation, coset graphs,
existence of Hamilton cycles in particular graphs, and group- and   finite field-theoretic results, respectively.

\subsection{Basic definitions and notation}
\label{ssec:definition}
\noindent
Throughout this paper, graphs are finite, simple and undirected,
and groups are finite.
Given a graph $\Gamma$, we denote the vertex set by $V(\Gamma)$ and the edge set by $E(\Gamma)$, respectively.
For two adjacent vertices $u,v \in V(\Gamma)$, we write
$u \sim v$ and denote the corresponding edge by $\{u,v\}$. For $\a\in V(\Gamma)$, we denote the neighborhood of $\a$ by $N(\a)$.
Let $U$ and $W$ be two disjoint subsets of $V(\Gamma)$. Then
by $X( U )$ and $X[U,W]$ we denote the subgraph of $\Gamma$ induced by $U$ and the bipartite subgraph of $\Gamma$ induced by the edges having one end-vertex in $U$
and the other end-vertex in $W$, respectively.

\vskip 3mm
Let $G$ be a group acting faithfully and transitively on a set $V$. A nonempty subset $B$ of $V$ is called a $block$ for $G$ if for each $g\in G$ either $B^g=B$ or $B^g\cap B=\emptyset$, where $V$ and the singletons $\{v\}$ $(v\in V)$ are called the {\it trivial} blocks, and other blocks are called {\it nontrivial} blocks. Put $\B=\{B^g: g\in G\}$. Then the sets in $\B$ form a partition of $V$. We call $\B$ the {\em system of blocks} containing $B$.
If $G$ has no nontrivial blocks on $V$, then $G$ is {\em primitive},
and is {\em imprimitive} otherwise. If every nontrivial normal subgroup of $G$ is transitive on $V$, then $G$ is {\em quasiprimitive}, and is {\em non-quasiprimitive} otherwise.
If $G$ acts non-quasiprimitively on $V$ with an intransitive normal subgroup $N$, then for any   $v\in V$, $v^N=\{v^g: g\in N\}$ is called an $N$-{\em block}. $|v^N|$ is called the length of blocks. The set of $N$-blocks forms a system of blocks for $G$.
\vskip 3mm
A graph $\Gamma$ is said to be {\it vertex-transitive} if its full automorphism group  $\Aut(\Gamma)$ acts transitively on $V(\Gamma)$. A vertex-transitive
graph $\Gamma$ is said to be {\it quasiprimitive} if there exists a group $G\le \Aut(\Gamma)$ such that $G$ acts quasiprimitively on $V(\Gamma)$ and is {\it non-quasiprimitive} otherwise. In particular, if there exists a group $G\le\Aut(\Gamma)$ such that $G$ acts primitively on $V(\Gamma)$, then $\Gamma$ is said to be  {\it primitive} and {\it imprimitive} otherwise.
\vskip 3mm
Let $\Gamma$ be a non-quasiprimitive graph and $G$ be a transitive subgroup of $\Aut(X)$. Then $G$ contains an intransitive normal subgroup $N$. Let $\B$ be the set of $N$-blocks and $A,B\in \B$. By $d(A)$ and $d(A,B)$ we denote the valency of $X(A)$ and $X[A,B]$, respectively. Since $N$ is transitive on both $V(X(A))$ and $V(X(B))$, we know that $X[A,B]$ is regular or is a null graph. A {\it normal quotient graph corresponding to} $N$ is the graph $\Gamma_\B$ whose vertex set is $\B$ with $A,B\in \B$ adjacent if there exist vertices $a \in A$ and $b \in B$, where $a,b\in V(\Gamma)$, such that $a \sim b$ in $\Gamma$. Let $N_0$ be the kernel of the action of $G$ on $\B$. $\Gamma_\B$ is also called the {\it $N$-block graph}. Then $N\le N_0$ and there exists a natural faithful action of $G/N_0$ on $\B$. In particular, the action is transitive. Thus, if $\Gamma$ is a connected non-quasiprimitive graph, then $\Gamma_\B$ is a connected vertex-transitive graph.

Let $(m,n)$ and $[m,n]$ denote the greatest common divisor and
the least common multiple of $m$ and $n$, respectively.

\subsection{Coset graphs and minimal connected coset graphs}
\label{ssec:orbital}
\subsubsection{Coset graphs}
Let $\Gamma$ be a connected vertex-transitive graph and  $G$ be a finite group acting transitively on $V(\Gamma)$, and fix an arbitrary vertex \(\alpha \in V(\Gamma)\). Write \(H = G_\alpha = \{g \in G : \alpha^g = \alpha\}\) for the stabilizer of \(\alpha\) in $G$. There is a natural $G$-equivariant bijection between the vertex set \(V(\Gamma)\) and the right coset space \([G:H] = \{Hg : g \in G\}\), which maps each vertex \(\alpha^g\) to the coset $Hg$ for all \(g \in G\). Under this identification, the transitive action of $G$ on vertices corresponds exactly to the right multiplication action of $G$ on cosets. Let $\Omega=\{g:\a^g\in N(\a)\}$ and $D=\cup_{g\in\Omega}HgH$. Then $\Gamma$ is determined by $G,H$ and $D$. This correspondence allows a purely group-theoretic construction of vertex-transitive graphs. The {\it coset graph} (or {\it double coset graph}) \(X(G, H, D)\) is defined to have vertex set \([G:H]\), with adjacency given by:
\(Hg_1 \sim Hg_2 \quad \text{if and only if} \quad Hg_2 g_1^{-1} \subseteq D.\) We denote the number of right cosets contained in $D$ by $|D|$, that is, $|D|=|\{Hg:g\in D\}|$. Similarly, $|Hg_0H|=|\{Hg:g\in Hg_0H\}|$ is called the length of the double coset $Hg_0H$.
\vskip 3mm
For $\Gamma=X(G,H,D)$ to be a finite, simple, undirected graph, three standard conditions hold: 
\begin{enumerate}
	\item \(H \cap D = \emptyset\) (no loops);
	\item  \(D = D^{-1} = \{d^{-1} : d \in D\}\) (symmetric adjacency for undirected edges);
	\item   $\Gamma$ is connected if $\lg g:g\in D\rg=G$.
\end{enumerate}
In particular, if $D_g=HgH=Hg^{-1}H$, then $D_g$ is said to be self-paired and is nonself-paired otherwise. If $D_g\neq D_{g^{-1}}$, then  $D_g$ and $D_{g^{-1}}$ are said to be symmetric-paired.

\vskip 3mm
$\Gamma=X(G,H,D)$ is $G$-symmetric provided $D=HgH=Hg^{-1}H$. $\Gamma=X(G,H,D)$ is $G$-edge transitive provided $D=HgH\cup Hg^{-1}H$, where $HgH\neq Hg^{-1}H$. 

\vskip 3mm

Let $\Gamma=X(G,H,D)$ be a connected vertex-transitive graph and $N$ be an intransitive normal subgroup of $G$. Suppose $|V(\Gamma)|=n$.
Then $N$ induces blocks of length $m$, where $m\ge 2$ and $m \mid n$. Let $\B$ be the set of $N$-blocks and $A, B \in \cal{B}$ be two distinct blocks. Let $N_0$ be the kernel of the action $G$ on $\B$. Then $\Gamma_\B\cong X(\og, \oh, \o{D})$, where $\overline{G}=G/N_0, \overline{H}=(HN_0)/N_0$ and $\overline{D}=\cup_{g\in D}\overline{HgH}$.

\subsubsection{Minimal connected coset graphs}

Let $\Gamma=X(G,H,D)$ be a coset graph, where $D=\cup_{g\in \Omega}HgH$.  Let $\Gamma_{0}=X(G,H,D_{0})$ be a subgraph of $\Gamma$ and $D_{g_{0}}=Hg_0H\cup Hg_0^{-1}H$, where $g_0\in D$ and $D_{0}=D\backslash D_{g_0}$. 	For any $g_0\in D$, if $\Gamma$ is connected but $\Gamma_{0}$ is not connected, then $\Gamma$ is called a {\it $G$-minimal connected (coset) graph} (in short, {\it $G$-minimal graph}). That is, if $\Gamma$ is a $G$-minimal graph, then $\Gamma_{0}=X(G,H,D_0)$ is not connected, where $D_0\subsetneq D$ and $D_0=D_0^{-1}$.

\vskip 3mm
Every connected vertex-transitive graph $\Gamma=X(G,H,D)$ contains at least one $G$-minimal graph. If every  $G$-minimal graph contains a Hamilton cycle, then every connected vertex-transitive graph with a transitive group $G$ contains a Hamilton cycle. Thus, to find a Hamilton cycle in $\Gamma$, we only need to consider $G$-minimal graphs.

\subsection{Existence of Hamilton cycles in particular graphs}
\label{ssec:numbers}
The following  known results about existence of Hamilton cycles in particular graphs will be used later.

\begin{proposition}\label{H-result}
	Let $p$ and $q$ be distinct primes, and $X$ a connected vertex-transitive graph of order $n$.
	Then one of the following holds.
	\begin{enumerate}
		\item If $G$ is a Cayley graph, then $X$ contains a Hamilton cycle when $n=p,pq$ or $2pq$ (\cite[Theorem~1.2]{MR234}).
		\item If $n=2p^2$, then $X$ contains a Hamilton cycle except for the Petersen graph (\cite{DM87}).
		\item If $n=pq$, then $X$ contains a Hamilton cycle except for the Petersen graph (\cite[Theorem~1.4]{MR4328721}).
		\item If $n=6p$ or $10p$, then $X$ contains a Hamilton cycle with the exception of the truncation of the Petersen graph (\cite[Theorem~1.1]{10p} and \cite[Theorem~1.1]{6p}).
		\item If $n=2pq$ and $X$ is not the Coxeter graph, then $X$ contains a Hamilton cycle when the automorphism group acts primitively on its vertices (\cite[Theorem~1.1]{MR4548744}).
	\end{enumerate}
\end{proposition}

\begin{proposition}
	\label{pro:2}
	{\rm \cite[Theorem~6]{BJ78}}
	Every $2$-connected regular graph of order $n$ and valency at least $n/3$ contains a Hamilton cycle.

\end{proposition}

\begin{proposition}
	\label{pro:3}
	{\rm \cite[Theorem~1.4]{MR4328721}}
	With the exception of the Petersen graph, a connected vertex-transitive graph of order $pq$, which $p$ and $q$ are primes, contains a Hamilton cycle.
\end{proposition}

\begin{proposition}
	\label{pro:4}
	{\rm \cite[Theorem~1.1]{10p},\cite[Theorem~1.1]{6p}}
	Every connected vertex-transitive graph of order $6p$ and $10p$, where $p$ is a prime, contains a Hamilton path. Moreover, with the exception of the truncation of the Petersen graph, every
	such graph contains a Hamilton cycle.
\end{proposition}

\begin{proposition}
	\label{pro:5}
	{\rm \cite[Theorem~1.1]{MR4548744}}
	Except for the Coxeter graph, every connected vertex-transitive graph of order $2rs$ contains a Hamilton cycle provided the automorphism group acts primitively on its vertices, where $r$ and $s$ are primes.
\end{proposition}

\subsection{Group- and  finite field-theoretic results}

Let $(m,n)$ and $[m,n]$ denote the greatest common divisor and
the least common multiple of $m$ and $n$, respectively.

Throughout, $\ZZ_n$ denotes a cyclic group of order $n$, and $D_{2n}$
denotes a dihedral group of order $2n$. We write finite field orders
as $r^n$, where $r$ is prime and $n\geq1$.
The notation $\ZZ_r^{\,a}$ denotes the direct product of $a$
copies of $\ZZ_r$; in particular,
$(\mathbb F_{r^n},+)\cong\ZZ_r^{\,n}$.
All parameters in the tables are subject to the stated primality
conditions, and all groups listed as nonabelian socles are assumed to
be nonabelian simple.
An action on cosets of $H\leq T$ means the action of $T$ on its right
cosets $H\backslash T$ by right multiplication.
Unless explicitly stated otherwise, the action columns describe the
restriction to $T=\operatorname{soc}(G)$.

For an odd prime $p$, every primitive permutation group of degree
$2p$ is $2$-transitive, except for $A_5$ and $S_5$ acting on the
ten unordered pairs of a five-element set; in particular, this
exception has $p=5$ and socle $A_5$ \cite{LS85}.
Consequently, apart from this degree-$10$ case, a simply primitive
group of degree $pq$, where $q<p$ are primes, has both $p$ and $q$
odd. We therefore list the odd-degree cases in Table~\ref{tab1}.

\begin{proposition}
	\label{pro:7}
	{\rm \cite{2TRAN}}
	Let $p$ and $q$ be odd primes with $q<p$, and let $G$ be a primitive
	permutation group of degree $pq$ that is not $2$-transitive.
	Then $G$ is almost simple, and its socle and the corresponding
	restricted action are as listed in Table~\ref{tab1}.
\end{proposition}

\begin{table}[htbp]
	\centering
	\caption{Socles and actions for simply primitive groups of degree
		$pq$, where $q<p$ are odd primes}
	\label{tab1}
	\renewcommand{\arraystretch}{1.18}
	\resizebox{\textwidth}{!}{%
		\begin{tabular}{|c|l|c|l|l|}
			\hline
			Row & $T=\operatorname{soc}(G)$ & $(p,q)$ & Action of $T$ & Conditions or remarks\\
			\hline
			1 & $A_p$ & $(p,(p-1)/2)$
			& unordered pairs of a $p$-set & $p\geq7$\\
			\hline
			2 & $A_{p+1}$ & $(p,(p+1)/2)$
			& unordered pairs of a $(p+1)$-set & $p\geq5$\\
			\hline
			3 & $A_7$ & $(7,5)$
			& unordered triples of a $7$-set & \\
			\hline
			4 & $\PSL(2,7)$ & $(7,3)$
			& cosets of $D_8$
			& $G=\operatorname{PGL}(2,7)$; $G_\alpha\cong D_{16}$\\
			\hline
			5 & $\PSL(4,2)$ & $(7,5)$
			& $2$-dimensional subspaces of $\mathbb F_2^4$ & \\
			\hline
			6 & $\PSL(5,2)$ & $(31,5)$
			& $2$- or $3$-dimensional subspaces of $\mathbb F_2^5$
			& dual actions\\
			\hline
			7 & $\operatorname{PSp}(4,2^{2^t})$
			& $(2^{2^{t+1}}+1,\,2^{2^t}+1)$
			& $1$-spaces or totally isotropic $2$-spaces
			& $t\geq1$\\
			\hline
			8 & $\operatorname{P\Omega}^{\varepsilon}(2d,2)$
			& $(2^d-\varepsilon,\,2^{d-1}+\varepsilon)$
			& singular $1$-spaces
			& $\varepsilon\in\{1,-1\}$, $d\geq4$\\
			\hline
			9 & $\PSL(2,r^2)$ & $((r^2+1)/2,r)$
			& cosets of $\operatorname{PGL}(2,r)$
			& $r$ an odd prime\\
			\hline
			10 & $\PSL(2,p)$ & $(p,(p-\varepsilon)/2)$
			& cosets of $D_{p+\varepsilon}$
			& $\varepsilon=\pm1$; $p\geq11$\\
			\hline
			11 & $\PSL(2,p)$
			& $(19,3),(29,7),(59,29),(61,31)$
			& cosets of $A_5$ & \\
			\hline
			12 & $\PSL(2,23)$ & $(23,11)$
			& cosets of $S_4$ & \\
			\hline
			13 & $\PSL(2,11)$ & $(11,5)$
			& cosets of $A_4$
			& $G=\operatorname{PGL}(2,11)$; $G_\alpha\cong S_4$\\
			\hline
			14 & $\PSL(2,13)$ & $(13,7)$
			& cosets of $A_4$ & \\
			\hline
			15 & $M_{11}$ & $(11,5)$
			& unordered pairs in its natural $11$-point action & \\
			\hline
			16 & $M_{22}$ & $(11,7)$
			& cosets of $\ZZ_2^4\rtimes A_6$ & \\
			\hline
			17 & $M_{23}$ & $(23,11)$
			& unordered pairs, or cosets of $\ZZ_2^4\rtimes A_7$
			& pairs in the natural $23$-point action\\
			\hline
	\end{tabular}}
\end{table}

In Rows~4 and~13 of Table~\ref{tab1}, the action of $G$ is primitive,
whereas its restriction to $T$ is imprimitive. In particular,
$D_{16}$ and $S_4$ in those rows are stabilizers in $G$, not in $T$.
The case $\operatorname{P\Omega}^{+}(6,2)\cong\PSL(4,2)$
is already covered by Row~5.

\begin{proposition}
	\label{pro:8}
	{\rm \cite{LS85}}
	Let $p$ and $q$ be primes with $q<p$, and let $G$ be a
	$2$-transitive permutation group of degree $pq$.
	Then $G$ is almost simple, and its socle and the corresponding
	action are as listed in Table~\ref{tab2}.
\end{proposition}

\begin{table}[htbp]
	\centering
	\caption{Socles and actions for $2$-transitive groups of degree
		$pq$, where $q<p$ are primes}
	\label{tab2}
	\label{tab3}
	\renewcommand{\arraystretch}{1.18}
	\resizebox{\textwidth}{!}{%
		\begin{tabular}{|c|l|c|l|l|}
			\hline
			Row & $T=\operatorname{soc}(G)$ & $(p,q)$ & Action of $T$ & Conditions or remarks\\
			\hline
			1 & $A_{pq}$ & $(p,q)$
			& natural action on $pq$ points
			& $A_{pq}\leq G\leq S_{pq}$\\
			\hline
			2 & $\PSL(d,r^n)$
			& $pq=\dfrac{r^{nd}-1}{r^n-1}$
			& projective points or hyperplanes
			& $r$ prime; $T\leq G\leq\operatorname{P\Gamma L}(d,r^n)$\\
			\hline
			3 & $\operatorname{PSU}(3,2^a)$
			& $(2^{2a}-2^a+1,\,2^a+1)$
			& isotropic $1$-spaces of $\mathbb F_{2^{2a}}^3$
			& $a\geq2$; both entries prime\\
			\hline
			4 & $\operatorname{Sz}(8)$ & $(13,5)$
			& the Suzuki ovoid of size $65$
			& $T\leq G\leq\operatorname{Aut}(T)$\\
			\hline
			5 & $A_7$ & $(5,3)$
			& cosets of $\PSL(3,2)$
			& $G=A_7$\\
			\hline
			6 & $M_{22}$ & $(11,2)$
			& cosets of $\PSL(3,4)$
			& natural $22$-point action; $T\leq G\leq\operatorname{Aut}(T)$\\
			\hline
	\end{tabular}}
\end{table}

In Row~2 of Table~\ref{tab2}, projective points are the
$1$-dimensional subspaces of $\mathbb F_{r^n}^d$, and hyperplanes are
its $(d-1)$-dimensional subspaces. When $q=2$, the degree condition
in that row forces $d=2$ and $r^n+1=2p$, with $r$ an odd prime;
the action is then on the projective line $\operatorname{PG}(1,r^n)$.
Thus the $2$-transitive groups of degree $2p$ are covered by
Rows~1, 2 and~6, with socles $A_{2p}$, $\PSL(2,r^n)$
and $M_{22}$, respectively.

\begin{proposition}
	\label{pro:9}
	{\rm \cite[Theorem~1.49]{fsg}}
	Let $G$ be a transitive permutation group of prime degree $p$.
	Then its socle $T$ and the corresponding restricted action
	are listed in Table~\ref{tab:prime-degree}.
	In Row~1, $G\leq\operatorname{AGL}(1,p)$ in its natural
	affine action on $\mathbb F_p$; in Rows~2--6, $G$ is almost
	simple with nonabelian simple socle $T$.
\end{proposition}

\begin{table}[htbp]
	\centering
	\caption{Socles and restricted actions for transitive groups
		of prime degree $p$}
	\label{tab:prime-degree}
	\renewcommand{\arraystretch}{1.18}
	\resizebox{\textwidth}{!}{%
		\begin{tabular}{|c|l|c|l|l|}
			\hline
			Row & $T=\operatorname{soc}(G)$ & $p$
			& Action of $T$ & Conditions or remarks\\
			\hline
			1 & $\ZZ_p$ & $p$
			& translations on $\mathbb F_p$
			& $G\cong\ZZ_p\rtimes\ZZ_e$, $e\mid p-1$\\
			\hline
			2 & $A_p$ & $p$
			& natural action on $p$ points & $p\geq5$\\
			\hline
			3 & $\PSL(2,11)$ & $11$
			& cosets of $A_5$ & \\
			\hline
			4 & $M_{11}$ & $11$
			& natural action on $11$ points & \\
			\hline
			5 & $M_{23}$ & $23$
			& natural action on $23$ points & \\
			\hline
			6 & $\PSL(d,r^n)$ & $\dfrac{r^{nd}-1}{r^n-1}$
			& projective points or hyperplanes
			& $d$ prime; $(d,r^n-1)=1$\\
			\hline
	\end{tabular}}
\end{table}

In Row~6 of Table~\ref{tab:prime-degree}, $r$ is prime and $n\geq1$.
For $z=(d,r^n-1)$, we have $p\equiv d\equiv0\pmod z$ and
$1\leq z\leq d<p$, so $z=1$ and $Z(\SL(d,r^n))=1$.
For projective degree $pq$, this only gives $z\mid pq$.

For Table~\ref{tab5}, a projective point stabilizer in
$\PSL(d,r^n)$ has the form
\[
P=\ZZ_r^{\,n(d-1)}\rtimes\bigl(\GL(d-1,r^n)/Z\bigr),
\qquad Z=Z(\SL(d,r^n)).
\]
Here the action is induced by
$A\mapsto\operatorname{diag}(\det(A)^{-1},A)$, and $Z$ is
identified with the scalar subgroup $\{\lambda I_{d-1}:\lambda^d=1\}$
of order $(d,r^n-1)$. Hyperplane stabilizers use the dual action.
For $L\leq\GL(d-1,r^n)/Z$ and
$H=\ZZ_r^{\,n(d-1)}\rtimes L$, we have
\[
[\PSL(d,r^n):P]=\frac{r^{nd}-1}{r^n-1},\qquad
[P:H]=[\GL(d-1,r^n)/Z:L].
\]

The next lemma adapts the prime-degree argument in the proof of
\cite[Lemma~4.5]{PWX93} to subgroups of index $2q$.

\begin{lemma}
	\label{lem:linear-index-2q}
	Let $d\geq3$ be an odd prime, $r$ a prime and $n\geq1$,
	with $p=(r^{nd}-1)/(r^n-1)$ prime.
	If $L\leq\GL(d-1,r^n)$ has index $2q$, where $p>q>2$
	are primes, then either $\SL(d-1,r^n)\leq L$ or
	\[
	(d,r^n,q)\in\{(3,2,3),(3,3,3),(3,5,3),(3,5,5)\}.
	\]
\end{lemma}

\demo Put $S=\SL(d-1,r^n)$ and suppose $S\not\leq L$.
For $(d-1,r^n)=(2,2)$ or $(2,3)$, the group orders force $q=3$.
Otherwise $S$ is perfect and $S/Z(S)$ is simple.
On the $2q$ cosets of $L$, the $S$-orbits have length $q$ or $2q$.
The central orbits have length at most $2$: larger ones would
give a nontrivial degree-$2$ action or an abelian transitive
image of $S$. Factoring them out gives $\PSL(d-1,r^n)$
a transitive action of degree $q$ or $2q$.

The prime-degree argument cited above, using
Proposition~\ref{pro:9} and primality of $p$, leaves only
$(d,r^n,q)=(3,5,5)$. An imprimitive degree-$2q$ action gives
the same conclusion by considering its $q$ blocks.
For a primitive action, Proposition~\ref{pro:8} and the
degree-$10$ exception give $d=3$, with $r^n\in\{4,5,9\}$
or $r^n+1=2q$. The first alternatives yield $p=21,31,91$,
so only $r^n=5$ and $q\in\{3,5\}$ survive.
In the latter, $q=3$ gives $r^n=5$; for $q>3$,
either $q\equiv1\pmod3$ and $3\mid p$, or $r=3$, $n$ is even,
and
\[
p=(3^n+3^{n/2}+1)(3^n-3^{n/2}+1),
\]
contrary to primality.\qqed

\begin{lemma}
	\label{max}
	Let $G$ be a primitive group of degree $n$, and let $M$ be a point stabilizer. Let $T$ be a maximal subgroup of $M$ and let $g\in G\backslash M$. If except for $T$, the group $M$ is the only nontrivial subgroup containing $T$, then $\lg T,g\rg =G$.
\end{lemma}

\demo Since $g\in G\backslash M$, we know $\lg T, g\rg \neq M$. Thus, $\lg T,g\rg =G$.\qqed

The argument of Praeger, Wang and Xu gives the following
orbit calculation.

\begin{proposition}
	\label{pro:linear-parabolic-orbits}
	{\rm \cite{PWX93}}
	Let $G=\PSL(d,r^n)$, where $d\geq3$, $r$ is prime and $n\geq1$,
	and let $M$ be a projective point stabilizer. Put
	$U=\ZZ_r^{\,n(d-1)}$ and suppose that
	\[
	U\rtimes\SL(d-1,r^n)\leq H<M<G
	\]
	in the standard parabolic embedding. In the action of $G$
	on $[G:H]$, the subgroup $H$ fixes the block
	$B=\{Hx:x\in M\}$ pointwise and is transitive on
	$[G:H]\setminus B$. Thus, putting $m=[M:H]$, the subdegrees are
	\[
	\underbrace{1,\ldots,1}_{m\ \mathrm{times}},\qquad |G:H|-m.
	\]
\end{proposition}

\demo Take $M$ to stabilize $\langle e_1\rangle$. As above,
\[
M\cong U\rtimes\bigl(\GL(d-1,r^n)/Z\bigr),
\qquad Z=Z(\SL(d,r^n)).
\]
The quotient of $M$ by $U\rtimes\SL(d-1,r^n)$ is cyclic
of order $(r^n-1)/(d,r^n-1)$. Hence $H\triangleleft M$ and
$m\mid(r^n-1)/(d,r^n-1)$.

Identify the vertices with nonzero vectors of $\FF_{r^n}^d$
modulo multiplication by elements of $(\FF_{r^n}^{\times})^m$.
The divisibility above ensures that this action factors through $G$.
The $m$ vertices represented by vectors in $\langle e_1\rangle$
form the block $B$, fixed pointwise by $H$.
The inverse image of $H$ in $\SL(d,r^n)$ contains
\[
\begin{pmatrix}1&u\\0&A\end{pmatrix},
\qquad
u\in\FF_{r^n}^{\,d-1},\quad A\in\SL(d-1,r^n).
\]
Since $d-1\geq2$, the matrices $A$ are transitive on nonzero
vectors of $\FF_{r^n}^{d-1}$, while $u$ adjusts the first
coordinate arbitrarily. Hence $H$ is transitive on all
vertices outside $B$, as required.\qqed

\begin{lemma}
	\label{re}
	Let $X=X(G,H,D)$ be a basic connected graph and  $\tau\in N_G(H)\backslash H$, where $\tau^q \in H$ and $q$ is a prime. Let $N=\lg H,\tau\rg$, $A=\{Hg:g\in Ng_1\}$ and $B=\{Hg:g\in Ng_2\}$, where $Ng_1\neq Ng_2$, where $g_2g_1^{-1}\in D$. Then $X[A,B]$ contains a matching.
\end{lemma}

\demo If $D=\tau^{-j}D\tau$, where $j\in \ZZ_q\backslash \{0\}$, then $Hg_1\sim Hg_2$ if and only if $H\tau^i g_1\sim H\tau^{ij} g_2$. Since $H\tau^{i_1j} g_2=H\tau^{i_2j} g_2$ if and only if $i_1=i_2$ and so we get a matching in $X[A,B]$.

\begin{lemma}
	\label{cong2}
	Let $\Gamma=X(G,H,D)$ be a connected coset graph and $G \triangleleft G_0$. Let $\tau\in G_0\backslash G$. Then $\Gamma\cong \Gamma_0= X(G,H^\tau,D^\tau)$.  In particular, let $N=N_{G_0}(H)$ and $l\in N\backslash H$. Then $\Gamma\cong X(G,H,D^l)$.
\end{lemma}

\demo Let $f$ be a bijection from $V(\Gamma)$ to $V(\Gamma_0)$, where $f(Hg)=H^\tau g^\tau$. Then $Hg_1\sim Hg_2$ if and only if $H^\tau g_1^\tau\sim H^\tau g_2^\tau$. 

\section{The proof of Theorem~\ref{th1.1}}


Let $G$ be a quasiprimitive imprimitive permutation group of degree
$2pq$, where $p>q>2$ are primes, and let $H_0$ be a point stabilizer.
Since $G$ is imprimitive, there exists a maximal proper subgroup
$M_0$ of $G$ such that $H_0<M_0<G$. As $[G:H_0]=2pq$, we have
\[
t=[G:M_0]\in\{2,q,p,pq,2q,2p\}.
\]
The action of $G$ on the right cosets of $M_0$ is primitive.
Moreover, its kernel $N=\operatorname{core}_G(M_0)$ is trivial.
Indeed, if $N\neq1$, then quasiprimitivity forces $N$ to be
transitive on the original permutation domain, giving
$G=H_0N\leq M_0$, a contradiction.
Thus $G$ admits a faithful primitive action of degree $t$.
Since $p\mid |G|$ and $G$ embeds in $S_t$, we have $t\geq p$,
and hence
\[
t\in\{p,pq,2q,2p\}.
\]

By the classification of primitive groups of prime degree or of
degree the product of two distinct primes, $G$ is either almost
simple or of affine type with $t=p$. The latter case is impossible,
since its normal subgroup of order $p$ cannot be transitive on
$2pq$ points. Therefore $G$ is almost simple, and
$T=\operatorname{soc}(G)$ is quasiprimitive and imprimitive on
the original permutation domain.

Set $H=H_0\cap T$ and $M=M_0\cap T$. Since $T$ is transitive,
$G=TH_0=TM_0$, and consequently
\[
H<M<T,\qquad
[T:H]=2pq,\qquad
[T:M]=[G:M_0]=t,\qquad
[M:H]=\frac{2pq}{t}.
\]

\begin{lemma}
	\label{le3.1}
	If $[T:M]$ is a product of two distinct primes, the coset
	action of $T$ on $[T:H]$ is among those described by the
	following subgroup chains.
	\begin{enumerate}
		\item $[T:M]=pq$ and $[M:H]=2$, with:
		\begin{enumerate}
			\item $T=A_p$, $M\cong S_{p-2}$, $H=A_{p-2}$,
			and $q=(p-1)/2$;
			\item $T=A_{p+1}$, $M\cong S_{p-1}$, $H=A_{p-1}$,
			and $q=(p+1)/2$;
			\item $T=\PSL(d,r^n)$, $pq=(r^{nd}-1)/(r^n-1)$,
			with $r$ odd and $d\geq3$ odd. Here
			\[
			M=\ZZ_r^{\,n(d-1)}\rtimes(\GL(d-1,r^n)/Z),\qquad
			H=\ZZ_r^{\,n(d-1)}\rtimes L,
			\]
			where $M$ stabilizes a projective point or hyperplane,
			$Z$ is as above, and $[\GL(d-1,r^n)/Z:L]=2$;
			\item $T=\PSL(2,r^2)$, $M\cong\operatorname{PGL}(2,r)$,
			$H\cong\PSL(2,r)$, and $(p,q)=((r^2+1)/2,r)$;
			\item $T=\PSL(2,p)$, $M\cong D_{p+\varepsilon}$,
			and $H\cong\ZZ_{(p+\varepsilon)/2}$ or
			$D_{(p+\varepsilon)/2}$, where $p\geq11$,
			$q=(p-\varepsilon)/2$ and $\varepsilon=\pm1$;
			\item $T=\PSL(2,23)$, $M\cong S_4$, $H\cong A_4$,
			and $(p,q)=(23,11)$;
			\item $T=M_{23}$, $M\cong\PSL(3,4).2$,
			$H\cong\PSL(3,4)$, and $(p,q)=(23,11)$,
			with $M$ stabilizing an unordered pair;
			\item $T=A_7$, $(p,q)=(7,5)$,
			$M=(S_3\times S_4)\cap A_7$, and
			$H=A_3\times A_4$, with $M$ stabilizing a $3$-subset;
			\item $T=\PSL(2,7)$, $(p,q)=(7,3)$,
			$M\cong D_8$, and $H\cong\ZZ_4$ or
			$D_4\cong\ZZ_2^2$;
			\item $T=\PSL(4,2)$, $(p,q)=(7,5)$,
			\[
			M=\ZZ_2^4\rtimes(S_3\times S_3),\qquad
			H=\ZZ_2^4\rtimes L,
			\]
			where $M$ stabilizes a $2$-space and
			$[S_3\times S_3:L]=2$;
			\item $T=\PSL(5,2)$, $(p,q)=(31,5)$,
			\[
			M=\ZZ_2^6\rtimes(S_3\times\PSL(3,2)),\qquad
			H=\ZZ_2^6\rtimes(\ZZ_3\times\PSL(3,2)),
			\]
			with $M$ stabilizing a $2$- or $3$-space;
			\item $T=M_{11}$, $(p,q)=(11,5)$,
			$M=\ZZ_3^2\rtimes R$, and $H=\ZZ_3^2\rtimes L$,
			where $R\in\operatorname{Syl}_2(\GL(2,3))$ and
			$[R:L]=2$. The three possibilities for $L$ are
			$\ZZ_8$, $D_8$ and the quaternion group $Q_8$;
			$M$ stabilizes an unordered pair.
		\end{enumerate}
		The pair stabilizers above refer to the natural actions.
		Both $T$-classes of $H$ are included in parts~(d) and~(f),
		as are both dihedral classes in parts~(e) and~(i).
		\item $[T:M]=2p$ and $[M:H]=q$, with:
		\begin{enumerate}
			\item $T=\PSL(2,r^n)$, $p=(r^n+1)/2$, $q\mid r^n-1$,
			\[
			M=\ZZ_r^{\,n}\rtimes\ZZ_{(r^n-1)/2},\qquad
			H=\ZZ_r^{\,n}\rtimes\ZZ_{(r^n-1)/(2q)};
			\]
			\item $T=A_5$, $(p,q)=(5,3)$, $M\cong S_3$,
			and $H\cong\ZZ_2$, with $M$ stabilizing an unordered pair.
		\end{enumerate}
		\item $[T:M]=2q$, $[M:H]=p$, and $p=2q-1$, with:
		\begin{enumerate}
			\item $T=A_{p+1}$, $M=A_p$, and $H=A_{p-1}$;
			\item $T=\PSL(2,p)$,
			$M=\ZZ_p\rtimes\ZZ_{(p-1)/2}$, and
			$H\cong\ZZ_{(p-1)/2}$.
		\end{enumerate}
	\end{enumerate}
	Different choices of the intermediate subgroup $M$ may
	yield the same coset action of $T$.
\end{lemma}

\demo Propositions~\ref{pro:7} and~\ref{pro:8} reduce $(T,M)$
to Tables~\ref{tab1} and~\ref{tab2}, together with the degree-$10$
exception with socle $A_5$. Each listed chain gives a faithful
quasiprimitive imprimitive action of $T$, since $T$ is simple.
We determine the subgroups of the required index in $M$.

\smallskip\noindent
\emph{(i) $t=pq$.}
An index-$2$ subgroup exists precisely when $|M/M'|$ is even, where $M'=[M,M]$ is the commutator subgroup of $M$.
The stabilizers in Rows~8, 11 and~16 of Table~\ref{tab1},
the subgroup $\ZZ_2^4\rtimes A_7$ in Row~17, and those in
Rows~1 and~5 of Table~\ref{tab2} are perfect.
Rows~7 and~13--14 of Table~\ref{tab1} have abelianizations
$\ZZ_{r^n-1}$ (with $r=2$) and $\ZZ_3$, respectively.
The unitary and Suzuki stabilizers in Rows~3--4 of
Table~\ref{tab2} have abelianizations
$\ZZ_{(r^{2n}-1)/(3,r^n+1)}$ and $\ZZ_7$, also of odd order.
Thus none admits an index-$2$ subgroup \cite{KL,Atlas}.
Row~6 of Table~\ref{tab2} has even degree.
For Row~2 of that table,
$M/M'\cong\ZZ_{(r^n-1)/(d,r^n-1)}$; hence $r$ and $d$
are odd, and the unique index-$2$ subgroup comes from the
Levi quotient $\GL(d-1,r^n)/Z$.

Rows~3, 5, 6 and~15 of Table~\ref{tab1} have abelianizations
$\ZZ_2$, $\ZZ_2^2$, $\ZZ_2$ and $\ZZ_2^2$, respectively.
Their index-$2$ subgroups are those listed: in the linear cases
and the $M_{11}$ case, take inverse images from the Levi factor
and $R$, respectively. In Row~4, $T=\PSL(2,7)$,
$M\cong D_8$ and $[M:H]=2$, so $H\cong\ZZ_4$ or
$D_4\cong\ZZ_2^2$. Both give faithful quasiprimitive
imprimitive actions of $T$ of degree $42=2\cdot7\cdot3$.
The actions with $H\cong D_4$ also occur in
Lemma~\ref{pro:prime-block-index}(5)(b), using $H<S_4<T$
with $[T:S_4]=7$.
For the dihedral family in part~(1)(e), write $M=\langle a,b\rangle$
with $|a|=(p+\varepsilon)/2$ and $b^2=1$. Its three index-$2$
subgroups are $\langle a\rangle$, $\langle a^2,b\rangle$ and
$\langle a^2,ab\rangle$; both dihedral subgroups are included.
The two dihedral subgroups above are not conjugate in $T$.
Part~(1)(d) includes both conjugacy classes in $T$ of
$H\cong\PSL(2,r)$, and part~(1)(f) includes both classes
of $H\cong A_4$.
The other cases in part~(1) follow from the index-$2$
subgroups of the remaining stabilizers $M$.

\smallskip\noindent
\emph{(ii) $t=2p$.}
The degree-$10$ exception gives $p=5$, $q=3$, $M\cong S_3$
and $H\cong\ZZ_2$, as in part~(2)(b).
Otherwise consider Rows~1, 2 and~6 of Table~\ref{tab2}.
Row~1 is excluded by simplicity and $|A_{2p-1}|>|S_q|$.
In Row~6, $p=11$ and $q\in\{3,5,7\}$, whereas
$\PSL(3,4)$ is simple of order $20160>|S_7|$.
For Row~2, put $U=\ZZ_r^n$, where $r^n=2p-1$.
If $U\not\leq H$, its image on the $q$ cosets is regular;
irreducibility of the square-scalar action makes this image
faithful, forcing $r^n=q$, a contradiction.
Thus $U\leq H$, and the cyclic quotient $M/U$ gives the
unique subgroup in part~(2)(a), requiring $q\mid r^n-1$.

\smallskip\noindent
\emph{(iii) $t=2q$.}
The degree-$10$ and $M_{22}$ cases are excluded because
$|T|$ has no prime divisor $p>q$.
For $M=A_{2q-1}$, Proposition~\ref{pro:9} forces $p=2q-1$.
For $T=\PSL(2,r^n)$ with $r^n+1=2q$, the conditions
$p>q$ and $p\mid|T|$ imply $p\mid r^n<2p$.
Thus $r=p$, $n=1$, and $H$ is a cyclic complement in
$\ZZ_p\rtimes\ZZ_{(p-1)/2}$, giving part~(3).\qqed
\begin{lemma}
	\label{pro:prime-block-index}
	If $[T:M]=p$, then $[M:H]=2q$ and one of the following holds.
	\begin{enumerate}
		\item $T=A_p$, $M=A_{p-1}$, with:
		\begin{enumerate}
			\item $p=2q+1$, $H\cong A_{p-2}$;
			\item $(p,q)=(5,3)$, $H\cong\ZZ_2$;
			\item $(p,q)=(7,5)$, $H\cong\ZZ_3^2\rtimes\ZZ_4$.
		\end{enumerate}
		Both classes of $A_5<A_6$ occur when $(p,q)=(7,3)$.
		\item $T=\PSL(2,11)$, $p=11$, $M\cong A_5$;
		$H\cong D_{10}$ for $q=3$, or $S_3$ for $q=5$.
		\item $T=M_{11}$, $(p,q)=(11,5)$,
		$M\cong M_{10}\cong A_6.2_3$, $H\cong\ZZ_3^2\rtimes Q_8$.
		\item $T=M_{23}$, $(p,q)=(23,11)$,
		$M\cong M_{22}$, $H\cong\PSL(3,4)$.
		\item $T=\PSL(d,r^n)$, $d\geq3$ prime, and
		\[
		p=\frac{r^{nd}-1}{r^n-1},\qquad
		M=\ZZ_r^{\,n(d-1)}\rtimes\GL(d-1,r^n),
		\]
		with $M$ a point or hyperplane stabilizer:
		\begin{enumerate}
			\item $H=\ZZ_r^{\,n(d-1)}\rtimes L$,
			$\SL(d-1,r^n)\leq L$, and $[\GL(d-1,r^n):L]=2q$;
			\item $T=\PSL(3,2)$, $(p,q)=(7,3)$,
			$M\cong S_4$, and $|H|=4$: $\ZZ_4$ or either class
			of $\ZZ_2^2$;
			\item $T=\PSL(3,3)$, $(p,q)=(13,3)$,
			$H=\ZZ_3^2\rtimes L$, with
			$L\cong\ZZ_8,D_8$ or $Q_8$ in $\GL(2,3)$;
			\item $T=\PSL(3,5)$, $p=31$, $H=\ZZ_5^2\rtimes L$.
			For $q=3$, $L$ stabilizes a $1$-space.
			For $q=5$, $L=\pi^{-1}(A_4)$ or $\pi^{-1}(S_3\times S_2)$,
			where $\pi:\GL(2,5)\to\operatorname{PGL}(2,5)\cong S_5$.
		\end{enumerate}
	\end{enumerate}
	Each chain gives a faithful quasiprimitive imprimitive action
	of $T$ of degree $2pq$.
	
\end{lemma}

\demo Proposition~\ref{pro:9} determines $(T,M)$.
Each listed chain gives the required action by its indices
and the simplicity of $T$. We exclude the other possibilities.

\smallskip\noindent
\emph{(i) Alternating groups.}
For $p=5$, necessarily $q=3$ and $|H|=2$.
For $p\geq7$, simplicity of $M=A_{p-1}$ and $|M|>q!$
exclude actions on $2$ or $q$ blocks, so its degree-$2q$
action is primitive. The classification gives $p-1=2q$,
$H\cong A_{p-2}$, or $M=A_6\cong\PSL(2,9)$ of degree $10$
with $H\cong\ZZ_3^2\rtimes\ZZ_4$.
Both index-$6$ classes of $A_5<A_6$ are included \cite{Atlas}.

\smallskip\noindent
\emph{(ii) Exceptional groups.}
In $M\cong A_5$, $q\in\{3,5\}$, with the required subgroups
$D_{10}$ and $S_3$, giving part~(2).
For $M\cong M_{10}$, $q\in\{3,5\}$ and the maximal subgroup
indices are $2,10,36,45$ \cite{Atlas}. An index-$6$ subgroup
would have index $3$ in $A_6$, which is impossible.
An index-$10$ subgroup is maximal, since $A_6$ has no
index-$5$ subgroup, and is isomorphic to $\ZZ_3^2\rtimes Q_8$.
For $M\cong M_{22}$, $q\in\{3,5,7,11\}$, while the least
proper-subgroup index is $22$, attained by $\PSL(3,4)$
\cite{Atlas}. Thus $q=11$.

\smallskip\noindent
\emph{(iii) Linear groups.}
Primality of $p$ gives $d$ prime, $(d,r^n-1)=1$, and
$M=U\rtimes\GL(d-1,r^n)$, where $U=\ZZ_r^{\,n(d-1)}$.
If $U\not\leq H$, its orbits have length $r\in\{2,q\}$,
and their point stabilizers are $\FF_r$-hyperplanes of $U$.
The linear factor permutes these hyperplanes transitively,
so counting $U$-orbits gives
\[
\frac{r^{n(d-1)}-1}{r-1}\ \bigm|\ \frac{2q}{r}.
\]
Thus $r=2$ and $q=2^{n(d-1)}-1$; primality of $p,q$
forces $(d,n)=(3,1)$ or $(2,2)$, giving $\PSL(3,2)$
with $|H|=4$ or $A_5$ with $|H|=2$.
If $d=2$ and $U\leq H$, then $M/U$ has odd order.
Hence we may assume $d\geq3$ odd, $H=U\rtimes L$ and
$[\GL(d-1,r^n):L]=2q$.

By Lemma~\ref{lem:linear-index-2q}, either $\SL(d-1,r^n)\leq L$,
giving part~(5)(a), or $d=3$ and $r^n\in\{2,3,5\}$.
For $\GL(2,2)$, $q=3$ and $L=1$.
For $\GL(2,3)$, $q=3$ and $|L|=8$; the three possibilities
$\ZZ_8,D_8,Q_8$ are the index-$2$ subgroups of a semidihedral
Sylow $2$-subgroup.

For $\GL(2,5)$, let $Z_0=Z(\GL(2,5))\cong\ZZ_4$.
If $Z_0\not\leq L$, then $[\GL(2,5):Z_0L]=3$ or $5$.
The former is impossible in $S_5$; in the latter,
$J=Z_0L$ maps onto $S_4$ and $[J:L]=2$.
The determinant maps $J$ onto $\mathbb F_5^\times$.
Also $J'$ maps onto $A_4$ and has order $24$:
order $12$ would embed $A_4$ in $\SL(2,5)$,
whose unique involution is $-I$.
Thus $J/J'\cong\ZZ_4$ via the determinant, so its unique
index-$2$ subgroup contains $Z_0$, a contradiction.
Hence $L=\pi^{-1}(\pi(L))$. For $q=3$, $\pi(L)$ has order $20$
and normalizes a Sylow $5$-subgroup; for $q=5$, it has
order $12$ and is conjugate to $A_4$ or $S_3\times S_2$.
This gives part~(5)(d).\qqed

\section{The proof of Theorem~\ref{the:main}}
We denote $\a\sim\b$ if $\a$ and $\b$ is adjoined.

\f {\bf Strategy of the proof of Theorem~\ref{the:main}:}  
To prove Theorem~\ref{the:main}, let $\Gamma=X(G,H,D)$ be a connected vertex-transitive graph of order $2pq$, where $G$ acts quasiprimitively on its vertices and $p,q$ are primes. 

By Proposition~\ref{H-result}, all connected vertex-transitive graphs of order $4p$, $2p^2$, $6p$ or $10p$ are Hamiltonian, except for the truncation of the Petersen graph. 
Consequently, we assume $p>q>5$ below.
If $G$ is primitive, then by Proposition~\ref{H-result} again, $\Gamma=X(G,H,D)$ is Hamiltonian, as desired.
Henceforth, we assume that $G$ is imprimitive.
Reasonably, we make the following hypothesis.

\begin{hypothesis}
	\label{hypothesis1}
	Let  $\Gamma=X(G,H,D)$ be a connected vertex-transitive graph of order $2pq$, where $p>q\ge7$ are primes,
	and $\Gamma$ also a $G$-minimal graph.
	Moreover, $G$ acts quasiprimitively but imprimitively on its vertices $V(\Gamma)$.
\end{hypothesis}

Since $G$ is a quasiprimitive group of degree $2pq$, we know that $\soc(G)$ is listed in Table~\ref{tab5}. 
We may assume $G=\soc(G)$ as $\soc(G)$ acts transitively on $\Gamma$. 
Next we verify the entries in Table~\ref{tab5} except for Rows 6–8,
see  Lemmas \ref{lemma4.1}, \ref{lemma4.1}, \ref{lemma4.1} and \ref{lemma4.1}, respectively.
Thus Theorem~\ref{the:main} is proved.
\subsection{$G\cong A_n$}
Firstly, we deal with the case $G\cong A_{p}$ or $A_{p+1}$, Rows 1 and 2 of Table~\ref{tab5}. 
In fact, we get a more general result.

\begin{theorem}\label{AnH-cycle}
	Let $H\leq G\cong A_n$ and $H\cong A_{n-2}$, where $n\ge 6$ is a positive integer.
	Then every connected graph $\Gamma=X(G,H,D)$ is Hamiltonian.
\end{theorem}
\demo
Suppose that  $\Gamma$ is $G$-minimal.
Let $\Omega=\{(x,y)\di x\neq y,~x,y=1,2\cdots,n\}$ be the set of order pairs.
Note that $G$ acts transitively on $\Omega$, and this action is equivalent to the action of $G$ on the coset $[G:H]$.
Let $V(\Gamma)=\Omega$ and $\a=(1,2)$ such that $G_{\a}=H$.

Next, we consider the neighborhood $N(\a)$ of $\a$. 
Suppose that $(3,4)\in N(\a)$. Then $\{(x,y)\di x\neq y,~x,y=3,\cdots,n\}=(3,4)^{H}\subseteqq N(\a)$.
Since $|(3,4)^{H}|=(n-2)(n-3)$, it follows $d(\Gamma)\ge(n-2)(n-3)\ge \frac{1}{3}|V(\Gamma)|$.
By Proposition~\ref{pro:2}, $\Gamma$ has a Hamilton cycle, as desired.
So in what follows, we assume that $(x,y)\notin N(\a)$, where $x,y=3,\cdots,n$ and $x\neq y$.

Suppose that $(3,1)\in N(\a)$.
Since the permutation $(123)\in G$ and $\{(3,1),(1,2)\}^{(123)}=\{(1,2),(2,3)\}$,
One yields $(2,3)\in N(\a)$.
Then $(3,1)^{H}\cup(2,3)^{H}=\{(x,1)\di x=3,\cdots,n\}\cup \{(2,y)\di y=3,\cdots,n\}\subseteqq N(\a)$,
which implies the edge set 
$E(\Gamma)\supseteq\{\{(i,j),(k,i)\},\{(j,i),(i,k)\}\di i\neq j,i\neq k,~{\rm and}~i,j,k=1,2,\cdots,n\}.$
Since $\Gamma$ is $G$-minimal and any subgraph of $\Gamma$ which contains the edge set $\{\{(i,j),(k,i)\},\{(j,i),(i,k)\}\di i\neq j,i\neq k,~{\rm and}~i,j,k=1,2,\cdots,n\}$ is connected, we get 
$$E(\Gamma)=\{\{(i,j),(k,i)\},\{(j,i),(i,k)\}\di i\neq j,i\neq k,~{\rm and}~i,j,k=1,2,\cdots,n\}.$$
Consider a serious of subsets $\Omega_i:=\{(i,k),(k,i)\di k=i+1,\cdots,n\}$, where $i=1,2,\cdots,n-3$.
Since the induced subgraph $X\lg \Omega_i\rg\cong K_{n-i,n-i}-K_2$, we get that for any two vertices $(i,k_1)$ and $(k_2,i)$,
there exists a Hamilton path $P_{(i,k_1,k_2)}$ in $X\lg \Omega_i\rg$.
Set two paths $P_1:(n-1,n)\sim(n-2,n-1)\sim(n,n-2)$ and $P_2:(n-2,n)\sim(n-1,n-2)\sim(n,n-1)$.
Then we get a Hamilton cycle $C_1$ of $\Gamma$, where
$$C_1:P_{(1,n,n)}, P_1, P_{(2,n,n)}, P_2, P_{(3,n,n)}, P_{(4,n,n)},\cdots, P_{(n-3,n,n)}, (1,n),$$
as desired.
So in what follows, we assume that $(x,1),(2,x)\notin N(\a)$, where $x=3,\cdots,n$.

Suppose that $(1,3)\in N(\a)$.
Then $\{(1,y)\di y=3,\cdots,n\}=(1,3)^{H}\subseteqq N(\a)$, which implies
$$E(\Gamma)\supseteq\{\{(i,j),(i,k)\}\di (i-k)(i-j)(j-k)\neq 0~{\rm and}~i,j,k=1,2,\cdots,n\}.$$
Reset $\Omega_i=\{(i,k)\di  k=1,\cdots,n~{\rm and}~k\neq i\}$, where  $i=1,\cdots,n$.
Since the induced subgraph $X\lg \Omega_i\rg\cong K_n$, we get that for any two vertices $(i,k_1)$ and $(i,k_2)$,
there exists a Hamilton path $P_{(i,k_1,k_2)}$ in $X\lg \Omega_i\rg$.
If $\{(1,y)\di y=3,\cdots,n\}=(1,3)^{H}=N(\a)$, then $\Gamma$ is not connected, a contradiction.
Hence, $(2,1)\in N(\a)$ or $(x,2)\in N(\a)$, where $x=3,\cdots,n$.
If $(2,1)\in N(\a)$, then 
$E(\Gamma)\supseteq\{ \{(i,j),(j,i)\}\di i\neq j~{\rm and}~i,j=1,2,\cdots,n\}.$
Thus, there exits a Hamilton cycle $C_2$ of $\Gamma$, where
$$C_2:P_{(1,2,n)}, P_{(n,1,n-1)}, P_{(n-1,n,n-2)}, P_{(n-2,n-1,n-3)},\cdots, P_{(3,4,2)}, P_{(2,3,1)},(1,2),$$
as desired.
Suppose that $(2,1)\notin N(\a)$ but $(3,2)\in N(\a)$.
Then $N(\a)=(1,3)^{H}\cup(3,2)^{H}=\{(1,y),(y,2)\di y=3,\cdots,n\}$,
which implies 
$E(\Gamma)=\{\{(i,j),(i,k)\},\{(i,j),(k,j)\}\di (i-k)(i-j)(j-k)\neq 0,~{\rm and}~i,j,k=1,2,\cdots,n\}.$
If $n$ is even, then there exits a Hamilton cycle $C_3$ of $\Gamma$, where
$$C_3:P_{(1,2,3)}, P_{(2,3,n)}, P_{(3,n,2)}, P_{(4,2,n)}, P_{(5,n,2)},\cdots
, P_{(n-2,2,n)}, P_{(n-1,n,3)}, P_{(n,3,2)},(1,2);$$
if $n$ is odd, then there exits a Hamilton cycle $C_4$ of $\Gamma$, where
$$C_3:P_{(1,2,3)}\sim P_{(2,3,n)}\sim P_{(3,n,2)}\sim P_{(4,2,n)}\sim P_{(5,n,2)}\sim\cdots
\sim P_{(n-2,2,n)}\sim P_{(n-1,n,3)}\sim P_{(n,3,2)}\sim(1,2),$$
as desired.
So in what follows, we assume that $(1,x)\notin N(\a)$, where $x=3,\cdots,n$.

Suppose that $(3,2)\in N(\a)$.
Then $\{(y,2)\di y=3,\cdots,n\}=(3,2)^{H}\subseteqq N(\a)$. 
If $N(\a)=\{(y,2)\di y=3,\cdots,n\}=(3,2)^{H}$, then $\Gamma$ is not connected.
Since $(1,x),(x,1),(2,x),(x,y)\notin N(\a)$, where $x,y=3,\cdots,n$ and $x\neq y$,
we get $N(\a)=\{(2,1),(y,2)\di y=3,\cdots,n\}$, which implies
$$E(\Gamma)=\{\{(i,j),(k,j)\},\{(i,j),(j,i)\}\di (i-k)(i-j)(j-k)\neq 0~{\rm and}~i,j,k=1,2,\cdots,n\}.$$
Reset $\Omega_i=\{(k,i)\di  k=1,\cdots,n~{\rm and}~k\neq i\}$, where  $i=1,\cdots,n$.
Since the induced subgraph $X\lg \Omega_i\rg\cong K_n$, we get that for any two vertices $(k_1,i)$ and $(k_2,i)$,
there exists a Hamilton path $P_{(i,k_1,k_2)}$ in $X\lg \Omega_i\rg$.
With the same argument as the case $(1,3)\in N(\a)$, there also exits a Hamilton cycle of $\Gamma$, as desired.
\qqed

By Theorem \ref{AnH-cycle}, the following result is obvious.

\begin{lemma}\label{lemma4.1}
	Under Hypothesis \ref{hypothesis1}, if $(G,H)=(A_p,A_{p-2})$ or $(A_{p+1},A_{p-1})$, that is Rows 1 and 2 of Table~\ref{tab5}, then $\Gamma$ is Hamiltonian.
\end{lemma}

\subsection{$G= \PSL(2,p^m)$}

Let $\lg \theta \rg=\FF_{p^m}=\FF_r$, the multiplicative group of the finite field $\FF_{r}$, where $r=p^m$. Set

$$\sigma_i=\begin{pmatrix}
	0 & \theta^i\\
	-\theta^{-i} & 0
\end{pmatrix},
t_i=\begin{pmatrix}
	1 & 0\\
	\theta^i & 1
\end{pmatrix},
u_i=\begin{pmatrix}
	1 & \theta^i\\
	0 & 1
\end{pmatrix},
\ell_i=\begin{pmatrix}
	\theta^i & 0\\
	0 & \theta^{-i}
\end{pmatrix},i\in \ZZ_{r-1}.
$$
Let $u=u_0$ and $t=t_0$.
One may assume that $H=\lg u_i,\ell_q\mid i\in \ZZ_m\rg$. Set $M=H\lg \ell_1 \rg$. Let $g=\begin{pmatrix}
	a&b\\
	c&d
\end{pmatrix}$. If $c=0$, then $Hg=H\begin{pmatrix}
	a&0\\
	0&a^{-1}
\end{pmatrix}=H\ell^a$. If $a=d=0$, then $Hg=H\ell^b t$. If $acd\neq 0$, then $Hg=\begin{pmatrix}
	d^{-1}&0\\
	c&d
\end{pmatrix}=H\ell^{-d}\begin{pmatrix}
	1&0\\
	\frac{c}{d}&1
\end{pmatrix}$.   

We denote $H\ell^i\begin{pmatrix}
	1&0\\
	x&1
\end{pmatrix}$ and $H\sigma^i P$ by $(\theta^i,x)$ and $(\theta^i,\infty)$, respectively.

In particular, if $q\mid \frac{p-1}{2}$, then there exists $x\in \FF_p$ such that $\begin{pmatrix}
	x&0\\
	0&x^{-1}
\end{pmatrix}\notin H$, we may assume that $\begin{pmatrix}
	k&0\\
	0&k^{-1}
\end{pmatrix}\notin H$; if $q\nmid \frac{p-1}{2}$, then for any $x\in \FF_p$, $\begin{pmatrix}
	x&0\\
	0&x^{-1}
\end{pmatrix}\notin H$.

\begin{lemma}
	\label{lemma4.2}
	Let $\Gamma=X(G,H,D)$ be a $G$-minimal graph. Then there exists $g\in G\backslash M$ such that $D=HgH\cup Hg^{-1}H$.
\end{lemma}
\demo Let $g_0\in G$. In particular, $M$ is the unique maximal subgroup of $G$ containing $H$.  If $g_0\in G\backslash M$, by Lemma~\ref{max}  then $\lg H,g_0 \rg=G$ and so $\Gamma=X(G,H,D_0)$ is a $G$-minimal graph, where $D_0=Hg_0H\cup Hg_0^{-1}H$. If $g_0\in M$, then $\lg H, g_0 \rg \leq M\neq G$ and so $X(G,H,D_M)$ is disconnected, where $D_M=\cup_{m\in M}HmH$. \qqed

\begin{lemma}
	\label{lemma4.3}
	$Hg\ell_1H=H\ell_1^{-1}gH$, where $g\in G\backslash M$.
\end{lemma}
\demo Set 

$$g=\begin{pmatrix}
	a & b\\
	c & d
\end{pmatrix}\in G\backslash M.$$
Clearly, $c\neq 0$. Then $Hg=Hg_1$, where 

$$g_1=\begin{pmatrix}
	0 & -c^{-1}\\
	c & d
\end{pmatrix}{\rm and}~ Hg_1H=\{Ha: a=\begin{pmatrix}
	0 & -c^{-1}\\
	c & x
\end{pmatrix}, x\in \FF_{r^m}\}.$$
Thus, $Hg_1\ell_1H=\{Ha: a=\begin{pmatrix}
	0 & (-c\theta)^{-1}\\
	c\theta & x
\end{pmatrix}, x\in \FF_{r^m}\}=H\ell_1^{-1}g_1H$.
\qqed

\begin{corollary}
	Let $\Gamma=X(G,H,D)$ be a $G$-minimal graph. Then $D=HgH$, where $g=\begin{pmatrix}
		0&\theta^{-c}\\
		\theta^{c}&0
	\end{pmatrix}$ and $c\in \ZZ_{{q}}\backslash\{0\}$. Let $\Omega =\{Hx:Hx\in D\}$. Then $|\Omega|=p$.
\end{corollary}

\begin{lemma}
	\label{conj1}
	Let $\Gamma_1=X(G,H,D_1)$ and $\Gamma_2=X(G,H,D_2)$, where $D_1=H\begin{pmatrix}
		0&1\\
		-1&0
	\end{pmatrix}H$ and $D_2=H\begin{pmatrix}
		0&c^{-1}\\
		-c&0
	\end{pmatrix}H$. Then $\Gamma_1\cong \Gamma_2$.
\end{lemma}

\demo Let $G_0\cong\PGL(2,r^m)$ and $G\triangleleft G_0$. Let $\tau=\begin{pmatrix}
	1&0\\
	0&\theta
\end{pmatrix}$. Then $\tau\in N_{G_0}(H)$ and $D_1^{\tau^i}=H\begin{pmatrix}
	0&\theta^i\\
	-\theta^{-i}&0
\end{pmatrix}$. By Lemma~\ref{cong2}, $\Gamma_1\cong \Gamma_2$. 

\begin{lemma}
	\label{lemma4.4}
	Let $X=X(G,H,D)$ be a basic connected graph. Then $X$ has a Hamilton cycle.
\end{lemma}
\demo By Lemma~\ref{conj1}, we only consider the case that $D=H\begin{pmatrix}.
	0&1\\
	-1&0
\end{pmatrix}H$. We first consider the neighbor of every vertex.

$$\begin{array}{lclclcl}
	N((\theta^{i},\theta^x))&=&\{Hg:g\in H\sigma H\ell^it_x\}&=&\{Hg:g\in H\ell^{-i}\sigma Ht_x\}\cup \{H\ell^{-i}\sigma t_x\}\\
	&=&\{Hg:g\in H\ell^{-i}\sigma u_yt_x\}\cup\{H\ell^{-i}\sigma t_x\}&=&\{Hg:g\in H\ell^{-i-y}t_{x}t_{-y}\}\cup\{H\ell^{-i}\sigma\}\\
	&=&\{(\theta^{-i-y},\theta^{-y}+\theta^x),(\theta^{-i},\infty):y\in \FF_{r^m}\}
	
\end{array}$$

Let $\olx =\{(i,x):i\in \ZZ_q\}$, where $x\in \Omega=\{0,1,\ldots, -1,\theta,\theta+1,\ldots,\theta^2,\theta^2+1\ldots,\theta^{m-1}-1,\infty\}$. Then by Lemma~\ref{re}, $X[\olx_1,\olx_2]\cong qK_{1,1}$ for any $x_1,x_2\in \Omega$. Let $\overline{X}$ be the quotient graph of $X$, where $V(\overline{X})=\{\olx:x\in \Omega\}$. Then $\overline{X}\cong K_{r^m+1}$. Let $\overline{C_1}$, $\overline{C_2}$ and $\overline{C_3}$ be 4 Hamilton cycles of $\overline{X}$, where $C_1$ is 

$$
\overline{1},\overline{2},\ldots, \overline{\theta^{m-1}-1},\overline{\infty},\overline{0},\overline{1},
$$
$C_2$ is 
$$
\overline{1},\overline{3},\overline{2},\overline{4},\overline{5}\ldots, \overline{\theta^{m-1}-1},\overline{\infty},\overline{0},\overline{1}
$$  

and $C_3$ is 

$$
\overline{1},\overline{k+1},\overline{k},\overline{k-1},\ldots,\overline{2},\overline{k+1},\overline{k+2},\ldots,\overline{\theta^m-1},\overline{\infty},\overline{0},\overline{1}
$$  

Suppose  $q\nmid p-1$. Then we lift $\overline{C_1}$ and $\overline{C_2}$ to a path of length $r+2$ or a cycle of length $r+1$, and denote these two as $C_1$ and $C_2$, respectively. We may asssume the start-vertex of $C_1$ and $C_2$ is $(1,1)$. Now, we consider the first 4 vertices of these two path. In $C_1$ is 
$$(1,1),(1,2),(1,3),(1,4)$$
and $C_2$ is 

$$(1,1),(2,3),(-\frac{1}{2},2),(-4,4).$$

Thus, one of these two paths must not be a cycle of length $r+1$. Suppose that $C_1$ is a path. Then $C_1$ is

$$
(1,1),\ldots,(\theta^i,1),
$$
where $i\neq 0$.
We use cosets to rewrite $C_1$, which is

$$
Hg_1,Hg_2,\ldots, Hg_r,Hg_{r+1},H\ell^ig_1.
$$

If the start-point is $H\ell^{ki}g_1$, where $k\in \ZZ_q$, then we can get the following paths which are

$$
H\ell^{ki}g_1,H\ell^{-ki}g_2,\ldots, H\ell^{ki}g_r,H\ell^{-ki}g_{r+1},H\ell^{(k+1)i}g_1.
$$

Let $C_{1,k}$ be the path which is

$$
H\ell^{ki}g_1,H\ell^{-ki}g_2,\ldots, H\ell^{ki}g_r,H\ell^{-ki}g_{r+1}.
$$
Then $$
C_{1,1},C_{1,2},\ldots,C_{1,q}
$$
is a Hamilton cycle of $X$.

Suppose $q\mid p-1$. Then using the same argument as the above, we know that one of the lifted paths of $C_1$ and $C_3$ must not be a cycle and so we can get a Hamilton cycle of $X$.\qqed

\subsection{$G= \PSL(d,r^m)$}

\begin{lemma}
	\label{lem:linear-parabolic-hamilton}
	Let $G=\PSL(d,r^n)$, where $d\geq3$, $r$ is prime and $n\geq1$,
	and let $M$ be a projective point stabilizer. Put
	$U=\ZZ_r^{\,n(d-1)}$. If
	\[
	U\rtimes\SL(d-1,r^n)\leq H<M<G,
	\] then the pair $(G,H)$
	determines a unique $G$-minimal graph $\Gamma=X(G,H,D)$, and this graph
	contains a Hamilton cycle.
\end{lemma}

\demo Take $M$ to stabilize $\langle e_1\rangle$, and put $m=[M:H]$.
As in the proof of Proposition~\ref{pro:linear-parabolic-orbits},
$H\triangleleft M$ and $M/H$ is cyclic of order $m$.
Let $\theta$ generate $\FF_{r^n}^{\times}$, and take
\[
\ell=\operatorname{diag}(\theta,\theta^{-1},1,\ldots,1),
\qquad
g=\begin{pmatrix}0&-1\\1&0\end{pmatrix}\oplus I_{d-2},
\]
where the matrices are identified with their images in $G$.
Then $g\notin M$, and $\ell H$ generates the cyclic group $M/H$.
Consequently,
\[
\langle H,\ell\rangle=M,\qquad
M=\bigsqcup_{i=0}^{m-1}H\ell^iH
=\bigsqcup_{i=0}^{m-1}H\ell^i.
\]
The natural projective action is $2$-transitive, so
$G=M\sqcup MgM$.
Since $d\geq3$, the element
\[
h=\operatorname{diag}(1,\theta,\theta^{-1},1,\ldots,1)
\]
belongs to $H$. Direct calculation gives
$g\ell=\ell^{-1}g$ and $\ell g=h^{-1}gh$, whence
\[
\ell^i g\ell^j=h^{\,j-i}g h^{\,i-j}\in HgH.
\]
It follows that $MgM=HgH$, and the double coset decomposition is
\[
G=\left(\bigsqcup_{i=0}^{m-1}H\ell^iH\right)\sqcup HgH.
\]
The cosets in $M$ form the block
$B=\{H\ell^i:0\leq i<m\}$, fixed pointwise by $H$.
Thus the subdegrees are
\[
\underbrace{1,\ldots,1}_{m\ \mathrm{times}},\qquad |G:H|-m.
\]

Every connected coset graph $X(G,H,D)$ must have $HgH\subseteq D$:
otherwise $D\subseteq M$ and $\langle H,D\rangle\leq M<G$,
so the graph is disconnected. Since $HgH=G\setminus M$,
the graph $\Gamma=X(G,H,HgH)$ is the complete multipartite graph
with $[G:M]$ parts of size $m$, and is connected.
Its edges form a single $G$-orbit, so it is $G$-minimal.
As every connected $G$-invariant graph on $[G:H]$ contains
$\Gamma$, this is the unique $G$-minimal graph for $(G,H)$.
Since $[G:M]=(r^{nd}-1)/(r^n-1)\geq3$, it is $2$-connected
and regular, with
\[
d(\Gamma)
=m([G:M]-1)
=\left(1-\frac{1}{[G:M]}\right)|G:H|
>\frac{|G:H|}{3}.
\]
By Proposition~\ref{pro:2}, $\Gamma$ contains
a Hamilton cycle.\qqed

The same conclusion holds for hyperplane stabilizers by duality.

\subsection{$G=M_{23}$}
$M_{23}$ has a natural action on $\Omega=\{(x,y):1\leq x\leq 23, 1\leq y\leq 23, x\neq y\}$, with the same suborbits as the action of $S_{23}$ on $\Omega$, by Theorem~\ref{AnH-cycle}, $\Gamma$ contains a Hamilton cycle.

\subsection{$G=\PSL(2,23)$}
Analysing all $G$-minimal graphs with Magma, we conclude that every \(\Gamma\) contains a Hamilton cycle.



{\footnotesize
	\begin{thebibliography}{99}
		\itemsep=0pt
		
		
		\bibitem{A81} B. Alspach,
		Lifting Hamilton cycles of quotient graphs, {\em Discrete Mathematics.} {\it 78} (1989), 25-36.
		
		\bibitem{metacy} B. Alspach,
		Hamilton Cycles in Metacirculant Graphs
		with Prime Power Cardinal Blocks, {\em Ann. Discrete Mathematics} {\bf 41} (1989), 7–16.
		
		\bibitem{A83} B. Alspach,
		The search for long paths and cycles in vertex-transitive graphs and digraphs, in:
		Kevin L. McAvaney (ed.), {\it Combinatorial Mathematics Vlll}, Lecture Notes in Mathematics {\bf 884},
		Springer-Verlag, Berlin, 1981, 14--22.
		
		
		\bibitem{bcc15}
		P.~J.~Cameron (ed.),
		Problems from the Nineteenth British Combinatorial Conference,
		{\em Discrete Math.} {\bf 167/168} (1997), 605--615.
		
		\bibitem{seven}
		P.~J.~Cameron, M.~Giudici, G.~A.~Jones, W.~M.~Kantor,
		M.~H.~Klin, D.~Maru\v{s}i\v{c}, L.~A.~Nowitz,
		Transitive permutation groups without semiregular subgroups,
		{\em J. London Math. Soc.} {\bf 66} (2002), 325--333.
		
		\bibitem{chen}
		C.~C.~Chen. On edge-Hamiltonian property of Cayley graphs, {\em Discrete Mathematics.} {\bf 72} (1988), 29-33.
		
		\bibitem{10p}
		H. Chen, J. J. Li and H. Yu, Hamilton Cycles In Vertex-Transitive Graphs of Order 10p, arXiv:2506.19888.
		
		\bibitem{C98}
		Y. Q. Chen,
		On Hamiltonicity of vertex-transitive graphs and digraphs of order $p\sp 4$,
		{\em J. Combin. Theory Ser. B} {\bf 72} (1998), 110--121.
		
		\bibitem{CG96}
		S. Curran and J. A. Gallian,
		Hamiltonian cycles and paths in Cayley graphs and digraphs - a survey,
		{\em Discrete Math.} {\bf 156} (1996), 1--18.
		
		\bibitem{DMMN07}
		E.~Dobson, A.~Malni\v c, D.~Maru\v si\v c and L.~A.~Nowitz,
		Semiregular automorphisms of vertex-transitive graphs of certain valencies,
		{\em J. Combin. Theory Ser. B} {\bf 97} (2007), 371--380.
		
		
		
		\bibitem{MR4328721}
		F. S. Du, K. Kutnar and D.~Maru\v si\v c,
		Resolving the Hamiltonian problem for vertex-transitive graphs of
		order a product of two primes,
		{\it Combinatorica} {\bf 41} (2021), 507--543.
		
		\bibitem{MR4548744}
		F. S. Du, Y. Tian and H. Yu,
		Hamilton cycles in primitive graphs of order $2rs$,
		{\it Ars Math. Contemp} {\bf 23} (2023), \# P3.05.
		
		\bibitem{6p}
		F. S. Du and T. L. Zhou, Hamilton cycles in vertex-transitive graphs of order $6p$, {\em Discrete Applied Mathematics.}  {\bf 368} (2025), 165–175.
		
		\bibitem{D83}
		E.~Durnberger,
		Connected Cayley graphs of semidirect products
		of cyclic groups of prime order by abelian groups are hamiltonian,
		{\em Discrete Math.} {\bf  46} (1983),  55--68.
		
		\bibitem{GWM14}
		E. Ghaderpour and D. W. Morris,
		Cayley graphs on nilpotent groups
		with cyclic commutator subgroup are hamiltonian,
		{\em Ars Math. Contemp.} {\bf 7} (2014), 55--72.
		
		\bibitem{G1}
		M.~Giudici,
		Quasiprimitive groups with no fixed point free elements of prime order,
		{\em J. London Math. Soc.} {\bf 67} (2003), 73--84.
		
		\bibitem{GM07}
		H. H. Glover and D. Maru\v si\v c,
		Hamiltonicity of cubic Cayley graph,
		{\em J. Eur. Math. Soc.} {\bf 9} (2007), 775--787.
		
		\bibitem{GKMM12}
		H. H. Glover, K. Kutnar, A. Malni\v c and D. Maru\v si\v c,
		Hamilton cycles in $(2,odd,3)$-Cayley graphs,
		{\em Proc. London Math. Soc.} {\bf 104} (2012), 1171--1197.
		
		\bibitem{fsg}
		D. Gorenstein,
		{\it Finite Simple Groups: An Introduction to Their Classifications,} Springer, New York, 1982.
		
		\bibitem{fgt}
		B. Huppert, {\it Endliche Gruppen. I,}
		Springer Berlin, Heidelberg, 1967.
		
		\bibitem{BJ78}
		B. Jackson,
		Hamilton cycles in regular 2-connected graphs.
		{\em J. Combin. Theory Ser. B.} {\bf 29}  (1980),   27--46.
		
		
		
		
		
		\bibitem{KM09}
		K.~Kutnar and D.~Maru\v si\v c,
		Hamilton cycles and paths in vertex-transitive graphs - Current directions,
		{\em Discrete Math.} {\bf 309} (2009), 5491--5500.
		
		\bibitem{KM08}
		K.~Kutnar and D.~Maru\v si\v c,
		Hamiltonicity of vertex-transitive graphs of order $4p$,
		{\em European J. Combin.}  {\bf 29} (2008), 423--438.
		
		
		\bibitem{KMZ12}
		K. Kutnar, D. Maru\v si\v c and C. Zhang,
		Hamilton paths in vertex-transitive graphs of order $10p$,
		{\em European J. Combin.} {\bf 33} (2012), 1043--1077.
		
		\bibitem{KS09}
		K. Kutnar and P. \v Sparl,
		Hamilton paths and cycles in vertex-transitive graphs of order $6p$,
		{\em Discrete Math.} {\bf 309} (2009), 5444--5460.
		
		
		\bibitem{LS85}
		M.~W.~Liebeck and J.~Saxl,
		Primitive permutation groups
		containing an element of large prime order,
		{\em J. London Math. Soc.} (2) {\bf 31} (1985), 237--249.
		
		\bibitem{L70}
		L.~Lov\' {a}sz,
		{ Combinatorial Structures and Their Applications,}
		ed. R.Guy, H.Hanam, N.Sauer and J.Schonheim, Gordon and Breach, New York, 1970.
		
		\bibitem{MR234}
		F. Maghsoudi,
		On Hamiltonicity of Cayley graphs of order $pqrs$,
		{\em Australas. J. Combin} {\bf 84} (2022), 124--166.
		
		\bibitem{DM83}
		D.~Maru\v si\v c,
		Hamiltonian circuits in Cayley graphs,
		{\em Discrete Math.} {\bf 46} (1983), 49--54.
		
		\bibitem{DM87}
		D.~Maru\v si\v c,
		Hamiltonian cycles in vertex symmetric graphs of order $2p\sp 2$,
		{\em Discrete Math.} {\bf 66} (1987), 169--174.
		
		\bibitem{DM92}
		D.~Maru\v si\v c,
		Hamiltonicity of vertex-transitive $pq$-graphs,
		{\em Fourth Czechoslovakian Symposium on Combin.,
			Graphs and Complexity} (1992), 209--212.
		
		\bibitem{M81}
		D.~Maru\v si\v c,
		On vertex symmetric digraphs,
		{\em Discrete Math.} {\bf 36} (1981), 69--81.
		
		
		\bibitem{DM85}
		D.~Maru\v si\v c,
		Vertex transitive graphs and digraphs of order $p\sp k$,
		Cycles in graphs (Burnaby, B.C., 1982) 115--128,
		{\em Ann. Discrete Math.} {\bf 27}, North-Holland, Amsterdam, 1985.
		
		\bibitem{MP82}
		D.~Maru\v si\v c and T.~D.~Parsons,
		Hamiltonian paths in vertex-symmetric graphs of order $5p$,
		{\em Discrete Math.} {\bf 42} (1982), 227--242.
		
		\bibitem{MP83}
		D.~Maru\v si\v c and T.~D.~Parsons,
		Hamiltonian paths in vertex-symmetric graphs of order $4p$,
		{\em Discrete Math.} {\bf 43} (1983), 91--96.
		
		\bibitem{MS4}
		D.~Maru\v si\v c and R.~Scapellato,
		A class of graphs arising from the action of $\PSL(2,q^2)$ on cosets  of $\PGL(2,q)$,
		{\em Discrete Math.} {\bf 134} (1994), 99--110.
		
		\bibitem{MS2}
		D.~Maru\v si\v c and R.~Scapellato,
		Characterizing vertex-transitive $pq$-graphs with an imprimitive subgroup of automorphisms,
		{\em J.~Graph Theory} {\bf 16} (1992), 375--387.
		
		\bibitem{MS5}
		D.~Maru\v si\v c and R.~Scapellato,
		Classifying vertex-transitive graphs whose order is a product of two primes,
		{\em Combinatorica} {\bf 14} (1994), 187--201.
		
		\bibitem{PWX93}
		C.~E.~Praeger, R.-J.~Wang and M.-Y.~Xu,
		Symmetric graphs of order a product of two distinct primes,
		{\em J. Combin. Theory Ser. B} {\bf 58} (1993), 299--318.
		
		\bibitem{2TRAN}
		K.~Meagher, P.~Spiga and P.~H.~Tiep,
		An Erdos--Ko--Rado theorem for finite $2$-transitive groups,
		{\em European J. Combin.} {\bf 55} (2016), 100--118.
		
		\bibitem{WM18}
		D. W. Morris,
		Cayley graphs on groups with commutator
		subgroup of order $2p$ are hamiltonian,
		{\em The Art of Discrete and Applied Math.} {\bf 1} (2018),
		$\#$P1.04.
		
		\bibitem{WM15}
		D. W. Morris,
		Odd-order Cayley graphs with commutator subgroup of order $pq$ are hamiltonian,
		{\em Ars Math. Contemp.} {\bf 8} (2015), 1--28.
		
		
		
		\bibitem{Z15}
		J.~Y.~Zhang,
		Vertex-transitive digraphs of order $p^5$ are Hamiltonian,
		{\em Electronic J. Combin.} {\bf 22} (2015), $\#$P1.76.
		
		\bibitem{Zh}
		T.~L.~Zhou,
		On Hamilton cycles in connected vertex-transitive graphs of order 
		$2pq$.
		
		\bibitem{p}
		Turner and James,
		Point-symmetric graphs with a prime number of points,
		{\em J. Combina. Theory} {\bf 3} (1967), 136-145.
		
		
		\bibitem{S12}
		J.~H.~Silverman, {\em A friendly introduction to number theory},
		4th Edition, Pearson, 2012.
		
		\bibitem{W66}
		H.~Wielandt,
		{\em Permutation groups},
		Academic Press, New York, 1966.
		
		\bibitem{KL}
		P. Kleidman and M. Liebeck,
		{\it The Subgroup Structure of the Finite Classical Groups,}
		Cambridge University Press, Cambridge, 1990.
		
		\bibitem{A89}
		B. Alspach,
		Lifting Hamilton cycles of quotient graphs,
		{\em Discrete Math.} {\bf 78} (1989), 25--36.
		
		\bibitem{AP82a}
		B.~Alspach and T.~D.~Parsons,
		On hamiltonian cycles in metacirculant graphs,
		{\em Ann. Discrete Math.} {\bf 15} (1982), 1--7.	
		
		\bibitem{MS3}
		D.~Maru\v si\v c and R.~Scapellato,
		A class of non-Cayley vertex-transitive graphs associated with $\PSL(2,p)$,
		{\em Discrete Math.} {\bf 109} (1992), 161--170.
		
		\bibitem{Atlas}
		J. H. Conway, R. T. Curtis, S. P. Norton, R. A. Parker and R. A. Wilson,
		\emph{Atlas of Finite Groups}, Clarendon Press, Oxford, 1985.
		
	\end{thebibliography}
}
\end{document}